\documentclass[12pt, twoside]{article}
\usepackage{amsmath,amsthm}
\usepackage{amsfonts}
\usepackage{amssymb,latexsym}
\usepackage{enumerate}
\newtheorem{theorem}{Theorem}[section]
\newtheorem{remark}[theorem]{Remark}
\newtheorem{proposition}[theorem]{Proposition}
\newtheorem{lemma}[theorem]{Lemma}
\newtheorem{definition}[theorem]{Definition}
\newtheorem{exam}[theorem]{Example}

\def\ind{1{\hskip -3 pt}\hbox{\textsc{I}}}
\def\n{\noindent}
\def\fr{\frac}
\def\Om{\Omega}
\def\E{\mathcal E}
\def\de{\delta}
\def\ve{\varepsilon}
\def\va{\varphi}

\def\F{\mathcal F}
\DeclareMathOperator{\loc}{loc}

\begin{document}
\setlength{\baselineskip}{18truept}
\pagestyle{myheadings}
\markboth{ N. V. Phu and N.Q. Dieu}{Maximal subextension and approximation of $m-$ subharmonic function}
\title {Maximal subextension and approximation of $m-$ subharmonic function}
\author{
 Nguyen Van Phu* and Nguyen Quang Dieu**
\\ *Faculty of Natural Sciences, Electric Power University,\\ Hanoi,Vietnam.\\
**Department of Mathematics, Hanoi National University of Education,\\ Hanoi, Vietnam;\\
Thang Long Institute of Mathematics and Applied Sciences,\\ 
Nghiem Xuan Yem, Hanoi, Vietnam
\\E-mail: phunv@epu.edu.vn\\ and ngquang.dieu@hnue.edu.vn}

\date{}
\maketitle

\renewcommand{\thefootnote}{}



\renewcommand{\thefootnote}{\arabic{footnote}}
\setcounter{footnote}{0}

\begin{abstract}
In this paper, we first study subextensions in the classes  $\mathcal{F}_{m}(\Omega)$ and $\mathcal{E}_{m,\chi}(\Omega)$. These results are then used to study approximation in the classes $\mathcal{F}_{m}(\Omega)$ and $\mathcal{E}_{m,\chi}(\Omega)$.
\end{abstract}

\section{Introduction}
Let $\Omega\subset \tilde{\Omega}$ be bounded domains in $\mathbb{C}^{n}$ and $u\in PSH(\Omega).$ A function $\tilde{u}\in PSH(\tilde{\Omega})$ is said to be a subextension of $u$ if $\tilde{u}\leq u$ on $\Omega.$ The subextension in the class $\mathcal{F}(\Omega),$ where $\Omega$ is a bounded hyperconvex domain in $\mathbb{C}^{n}$ has been studied by Cegrell and Zeriahi in $\cite{CZ03}.$ In details, the authors proved that if $\Omega\Subset \tilde{\Omega}$ are bounded hyperconvex domains in $\mathbb{C}^{n}$ and $u\in\mathcal{F}(\Omega)$ then there exists $\tilde{u}\in\mathcal{F}(\tilde{\Omega})$ such that $\tilde{u}\leq u$ on $\Omega$ and $\int_{\tilde{\Omega}}(dd^{c}\tilde{u})^{n}\leq \int_{\Omega}(dd^{c}u)^{n}.$ For the class $\mathcal{E}_{p}(\Omega),p>0,$ the subextension was studied by Hiep. In $\cite{H08}$, it was shown that if $\Omega\subset \tilde{\Omega}$ are  hyperconvex domains in $\mathbb{C}^{n}$ and $u\in\mathcal{E}_{p}(\Omega)$, then there exists $\tilde{u}\in\mathcal{E}(\tilde{\Omega})$ such that $ \tilde{u}\leq u$  on $\Omega$ and $\int_{\tilde{\Omega}}(-\tilde{u})^{p}(dd^{c}\tilde{u})^{n}\leq \int_{\Omega}(-u)^{p}(dd^{c}u)^{n}.$ Using the subextension in the class $\mathcal{F}(\Omega),$ in $\cite{CH08},$ Cegrell and Hed studied the approximation in $\mathcal{F}(\Omega)$ by function in $\mathcal{F}(\Omega_j)$ where $\Omega, \Omega_j$ are hyperconvex domains such that $\Omega\Subset\Omega_{j+1}\Subset\Omega_j$ for all $j.$ More precisely, the authors proved the following interesting result: Let $u\in\mathcal{F}_{m}(\Omega)$ and $u_j$ be the maximal subextention of $u$ to $\Omega_j.$  Then $u_j$ increases to $u$ q.e on $\Om$ is true for all $u\in\mathcal{F}(\Omega)$ if it is true for one function $u\in\mathcal{F}(\Omega).$
Next, in $\cite{Be11},$ Benekourchi proved that any function in a class $\mathcal{E}_{\chi}$  can be approximated from below by an increasing sequence of plurisubharmonic functions defined on larger domains having finite $\chi-$ energy.
   
Motivated by the above results, we would like to carry out the same framework where the class of plurisubharmonic
functions is replaced by $m-$subharmonic ones. The paper is organized as follows. Besides the introduction,  the paper has other four sections. In Section 2 we recall the definitions and results concerning the $m$-subharmonic functions which were introduced and investigated intensively in the recent years by many authors (see \cite{Bl1}, \cite{SA12}). We also recall the Cegrell classes of $m$-subharmonic functions $\mathcal{F}_m(\Omega)$, $\mathcal{N}_{m}(\Omega),$  $\mathcal{E}_m(\Omega)$ and $\mathcal{E}_{m,\chi}(\Omega)$ which were introduced and studied in $\cite{Ch12}$, $\cite{T18}, \cite{DBH14}$ and $\cite{DT23}$.  In Section 3, we present subextension in the class $\mathcal{E}_{m}^{0}(\Omega),\mathcal{F}_{m}(\Omega)$ and $\mathcal{E}_{m,\chi}(\Omega)$. 
A key result is Lemma \ref{lm3.2} about  maximal subextension of a function $u \in \mathcal{E}_{m}^{0}(\Omega)$.
It is worth remarking that $u$ is not assumed to be continuous as in the plurisubharmonic case (Lemma 4.4 in $\cite{H08}$). From this lemma, we obtain maximal subextensions for elements in  $\mathcal{E}_{m,\chi}(\Omega).$
One novelty of our result is we do not require relative compactness of the smaller domain inside the larger one.
Finally, in Section 4, we use the obtained results  to study approximation of functions in $\mathcal{F}_{m}(\Omega)$ and $\mathcal{E}_{m,\chi}(\Omega)$ by functions in Cegrell classes defined on larger domains. The first main result of this section is
Theorem 4.1 where we generalize the mentioned above theorem of Cegrell and Hed  
in $\cite{CH08}$ to the class $\mathcal{F}_{m}(\Omega)$. It is worth mentioning that in Theorem 4.1, at some point, we allow the sequence $\Omega_j$ to be {\it non-decreasing}, at a price, we can only conclude the approximation occurs in the sense of convergence in capacity (see the implication $(iv) \Rightarrow (i)$).
The other main result is Theorem 4.4 where we study approximation of a function in the weighted energy class $\mathcal{E}_{m,\chi}(\Omega)$ by an increasing sequence of $m-$subharmonic functions 
in $\mathcal{E}_{m,\chi}(\Omega_j)$
where $\Omega_j \supset \Omega$ is a decreasing sequence of $m-$hyperconvex domains.\\ 
\section{Preliminaries}
  In this section, we first recall basic properties of $m$-subharmonic functions introduced by B{\l}ocki in \cite{Bl1} with the emphasis on the definition of the Hessian operator. This development is analogous to basic material for plurisubharmonic functions that can be found in \cite{BT1}, \cite{Kl},... Later, we repeat the construction 
of the Cegrell classes $\mathcal{E}^0_m(\Omega)$, $\mathcal{F}_m(\Omega)$  introduced and investigated by  Chinh more recently in \cite{Ch12} (see also \cite{SA12}). It should be noticed that many of the results for Cegrell classes of $m-$subharmonic functions are similar to those studied much earlier by Cegrell in \cite{Ce98}, \cite{Ce04} in the context of plurisubharmonic functions.

\n 
Let $\Omega$ be an open subset in $\mathbb{C}^n$.  By $\beta= dd^c\|z\|^2$ we denote the canonical K\"ahler form of $\mathbb{C}^n$ with the volume element $dV_{2n}= \frac{1}{n!}\beta^n$ where $d= \partial +\overline{\partial}$ and $d^c =\frac{\partial - \overline{\partial}}{4i}$.

\noindent{\bf 2.1} First, we recall the class of $m$-subharmonic functions which were introduced and investigated in \cite{Bl1}. For $1\leq m\leq n$, we define
$$\widehat{\Gamma}_m=\{ \eta\in {\mathbb{C}}_{(1,1)}: \eta\wedge \beta^{n-1}\geq 0,\ldots, \eta^m\wedge \beta^{n-m}\geq 0\},$$
where ${\mathbb{C}_{(1,1)}}$ denotes the space of $(1,1)$-forms with constant coefficients.

\noindent\begin{definition}\label{dn1}{\rm Let $u$ be a subharmonic function on an open subset $\Omega\subset \mathbb{C}^n$. Then $u$ is said to be an {\it $m$-subharmonic} function on $\Omega$ if for every $\eta_1,\ldots, \eta_{m-1}$ in $\widehat{\Gamma}_m,$ we have the inequality
$$ dd^c u \wedge\eta_1\wedge\cdots\wedge\eta_{m-1}\wedge\beta^{n-m}\geq 0$$
holds in the sense of currents.}
\end{definition}

\noindent   We denote $SH_m(\Omega)$ the set of $m$-subharmonic functions on $\Omega$ and $SH^{-}_m(\Omega)$ denotes the set of negative $m$-subharmonic functions on $\Omega$. It is clear that if $u\in SH_m$ then $dd^cu\in \widehat{\Gamma}_m$.

\n Now assume that $\Omega$ is an open set in $\mathbb{C}^n$ and $u\in \mathcal{C}^2(\Omega)$. Then from the Proposition 3.1 in \cite{Bl1} (also see the Definition 1.2 in \cite{SA12}) we note that $u$ is $m$-subharmonic function on $\Omega$ if and only if $(dd^c u)^k\wedge\beta^{n-k}\geq 0,$ for $k=1,\ldots, m$. More generally, if $u_1, \ldots, u_k\in\mathcal{C}^2(\Omega),$ then for all $\eta_1, \ldots, \eta_{m-k}\in \widehat{\Gamma}_m$, we have

\begin{equation*}\label{pt1}
dd^cu_1\wedge\cdots\wedge dd^cu_k\wedge\eta_1\wedge\cdots\wedge\eta_{m-k}\wedge\beta^{n-m}\geq 0
\end{equation*}
holds in the sense of currents. 

\n 
Note that, by $\cite{Ga59}$, a smooth function $u$ defined in open subset $\Omega\in\mathbb{C}^{n}$ is $m-$ subharmonic if 
$$S_{k}(u)=\sum_{1\leq j_{1}<...<j_{k}\leq n}\lambda_{j_{1}}\cdots\lambda_{j_{k}}\geq 0, \forall k=1,\cdots,m$$
where $(\lambda_{j_{1}}\cdots\lambda_{j_{k}})$ is the eigenvalue vector of the complex Hessian metrix of $u.$
\begin{exam}
Let $u(z)=4|z_{1}|^{2} +5|z_{2}|^{2}-|z_{3}|^{2}.$ Then we have $u\in SH_{2}(\mathbb{C}^{3}).$ However, $u$
 is not a plurisubharmonic function in $\mathbb{C}^{3}$ because the restriction of $u$
 on the line $(0,0,z_{3})$ is not subharmonic.
\end{exam}
We collect below basic properties of $m$-subharmonic functions that might be deduced directly from Definition 2.1. For more details, the reader may consult  $\cite{Ch15}, \cite{DBH14}, \cite{SA12}.$

\begin{proposition} \label{basic}
	 Let $\Om$ be an open set in $\mathbb C^n$. Then the following assertions hold true:

\n 	
(1) If $u,v\in SH_{m}(\Omega)$ then $au+bv \in SH_{m}(\Omega) $ for any $a,b\geq 0.$
\newline
(2) $PSH(\Omega)= SH_{n}(\Omega)\subset \cdots\subset SH_{1}(\Omega)=SH(\Omega).$
\newline
(3) If $u\in SH_{m}(\Omega)$ then a standard approximation convolution $ u*\rho_{\varepsilon}$ is also an m-subharmonic function on $\Omega_{\varepsilon}=\{z\in\Omega: d(z,\partial\Omega)>\varepsilon\}$ and $u*\rho_{\varepsilon}\searrow u$  as $\varepsilon\to 0.$ 
\newline
(4) The limit of an uniformly converging or decreasing sequence of $m$-subharmonic function is  $m$-subharmonic. 
\newline
(5) Maximum of a finite number of $m$-subharmomic functions is a $m$-subharmonic function.
\end{proposition}
\noindent Now as in \cite{Bl1} and \cite{SA12} we define the complex Hessian operator for locally bounded $m$-subharmonic functions as follows.
\begin{definition}\label{dn2}{\rm Assume that $u_1,\ldots, u_p\in SH_m(\Omega)\cap L^{\infty}_{\loc}(\Omega)$. Then the {\it complex Hessian operator} $H_m(u_1,\ldots,u_p)$ is defined inductively by
$$dd^cu_p\wedge\cdots\wedge dd^cu_1\wedge\beta^{n-m}= dd^c(u_p dd^cu_{p-1}\wedge\cdots\wedge dd^cu_1\wedge\beta^{n-m}).
$$}
\end{definition}
\noindent It was shown in \cite{Bl1} and later in \cite{SA12} that $H_m(u_1,\ldots, u_p)$ is a closed positive current of bidegree $(n-m+p,n-m+p).$ Moreover, this operator is continuous under decreasing sequences of locally bounded $m$-subharmonic functions. In particular, when $u=u_1=\cdots=u_m\in SH_m(\Omega)\cap L^{\infty}_{\loc}(\Omega)$ the Borel measure $H_m(u) = (dd^c u)^m\wedge\beta^{n-m}$ is well defined and is called the complex $m$-Hessian of $u$.

\noindent{\bf 2.2} Next, we recall the classes $\mathcal{E}^0_m(\Omega)$, $\mathcal{F}_m(\Omega)$ and $\mathcal{E}_m(\Omega)$ introduced and investigated in \cite{Ch12}. Let $\Omega$ be a bounded $m$-hyperconvex domain in $\mathbb{C}^n$, which mean there exists an $m-$ subharmonic function $\rho:\Omega\to (-\infty,0)$ such that the closure of the set $\{z\in\Omega:\rho(z)<c\}$ is compact in $\Omega$ for every $c\in (-\infty,0).$ Such a function $\rho$ is called the exhaustion function on $\Omega.$ If moreover function $\rho$ is defined in a neighborhood $\Omega^{'}$ of $\overline{\Omega}$ and $\Omega=\{\rho<0\}$ then we say that $\Omega$ is strongly $m-$hyperconvex domain.
By a recent result in \cite{AL}, for every bounded $m-$hyperconvex domain $\Om,$ there exists a negative, smooth, strictly $m-$subharmonic exhaustion function $\rho$ on $\Om.$ Here the strict $m-$subharmonicity means that for every point $a \in \Om$ there is a constant $c_a>0$ such that $\va-c_a |z|^2$ is $m-$subharmonic on a neighbourhood of $a.$
This implies that $\Omega$ admits an exhaustion of relatively compact strictly $m-$pseudoconvex domains $\{\Om_j\}$.

Throughout this paper, unless otherwise specified, by $\Omega$ we always mean a bounded  $m-$ hyperconvex domain in $\mathbb{C}^{n}.$
 Put
$$
\mathcal{E}^0_m=\mathcal{E}^0_m(\Omega)=\{u\in SH^{-}_m(\Omega)\cap{L}^{\infty}(\Omega):
\underset{z\to\partial{\Omega}}
\lim u(z)=0, \int\limits_{\Omega}H_m(u) <\infty\},$$
$$\mathcal{F}_m=\mathcal{F}_m(\Omega)=\big\{u\in SH^{-}_m(\Omega): \exists \mathcal{E}^0_m\ni u_j\searrow u, \underset{j}\sup\int\limits_{\Omega}H_m(u_j)<\infty\big\},$$
and
\begin{align*}
&\mathcal{E}_m=\mathcal{E}_m(\Omega)=\big\{u\in SH^{-}_m(\Omega):\forall z_0\in\Omega, \exists \text{ a neighborhood } \omega\ni z_0,  \text{ and }\\
&\hskip4cm \mathcal{E}^0_m\ni u_j\searrow u \text{ on } \omega, \underset{j}\sup\int\limits_{\Omega}H_m(u_j) <\infty\big\}.
\end{align*}
\n In the case $m=n$ the classes $\mathcal{E}^0_m(\Omega)$, $\mathcal{F}_m(\Omega)$ and $\mathcal{E}_m(\Omega)$ coincide, respectively, with the classes $\mathcal{E}^0(\Omega)$, $\mathcal{F}(\Omega)$ and $\mathcal{E}(\Omega)$ introduced and investigated earlier by Cegrell in \cite{Ce98} and \cite{Ce04}.

\noindent From  Theorem 3.14 in \cite{Ch12} it follows that if $u\in \mathcal{E}_m(\Omega)$, the complex $m$-Hessian $H_m(u)= (dd^c u)^m\wedge\beta^{n-m}$ is well defined and it is a Radon measure on $\Omega$.
Moreover, according to Theorem 3.18 in \cite{HP17}, the complex $m-$ Hessian operator is continuous under monotone convergence sequences in $\E_m$. Slightly weaker versions of this useful result for locally bounded $m-$subharmonic functions were proved earlier in  Theorem 2.6 and Theorem 2.10 of \cite{Ch15}.

\begin{exam}{\rm  For $0 < \alpha < 1$ we define the function
$$u_{m,\alpha} (z):= - (-\log \|z\|)^{\frac{\alpha m}{n}} + (\log 2)^{\frac{\alpha m}{n}}, 1\leq m\leq n, $$
on the ball $\Omega:= \{z\in\mathbb C^{n}:\|z\|<\frac 1 2\}$. 
Direct computations as in Example 2.3 of \cite{Ce98} shows that 
$u_{m,\alpha}\in\mathcal E_{m}(\Omega)$, $\forall 0<\alpha < \frac 1 m$. 

}\end{exam}
\noindent{\bf 2.3.} 
We say that an $m-$ subharmonic function $u$ is $m-$maximal if for every relatively compact open set $K$ on $\Omega$ and for each upper semicontinuous function  $v$ on $\overline{K},$ $v\in SH_{m}(K)$ and $v\leq u$ on $\partial K,$ we have $v\leq u$ on $K.$ The family of $m-$maximal $m-$ subharmonic function defined on $\Omega$ will be denoted by $MSH_{m}(\Omega).$ 

\n {\bf 2.4.} Following \cite{Ch15}, a set $E\subset\mathbb{C}^n$ is called $m$-polar if  $E\subset \{v=-\infty\}$ for some $v\in SH_m(\mathbb{C}^n)$ and $v \not\equiv -\infty.$
A property $P$ is said to be true quasi everywhere (q.e.) on a subset $E$ of $\mathbb C^n$ if $P$ holds true on $E \setminus F$ where $F$
is a $m-$polar subset of $\mathbb C^n.$

\noindent{\bf 2.5.} In the same fashion as the relative capacity introduced by Bedford and Taylor in \cite{BT1}, the $Cap_m$ relative capacity is defined as follows.
\begin{definition}\label{dn61}{\rm Let $E\subset \Omega$ be a Borel subset. The $m$-capacity of $E$ with respect to $\Omega$ is defined in \cite{Ch15} by
		$$Cap_{m}(E)=Cap_m(E,\Omega)= \sup\Bigl\{\int\limits_{E}H_m(u): u\in SH_m(\Omega), -1\leq u\leq 0\Bigl\}.$$}
\end{definition}

\noindent Proposition 2.8 in \cite{Ch15} gives some elementary properties of the $m$-capacity similar to those presented in \cite{BT1}. Namely, we have:

a) $Cap_m(\bigcup\limits_{j=1}^\infty E_j)\leq \sum\limits_{j=1}^\infty Cap_m(E_j)$.

b) If $E_j\nearrow E$ then $Cap_m(E_j)\nearrow Cap_m(E)$.\\
As in $\cite{BT1},$ by Theorem 2.20 in $\cite{Ch15}$ we have that 
$$Cap_{m}(E)=\int_{\Omega}H_{m}(h^{*}_{m,E,\Omega})$$
where $(h^{*}_{E,\Omega})$  is the smallest upper semicontinous majorant of $h_{E,\Omega}$ with
$$h_{m,E,\Omega}(z)=\sup\{\varphi(z): \varphi\in SH^{-}_{m}(\Omega):\varphi\leq -1 \,\,\text{on}\,\, E\}.$$\\
Note that, by Proposition 2.17 in $\cite{Ch15}$ if $\Omega$ be a $m-$ hyperconvex domain in $\mathbb{C}^{n}$ and $E\Subset\Omega$ then $h^{*}_{m,E,\Omega}\in \mathcal{E}_{m}^{0}(\Omega)$ and by Theorem 5.1 in $\cite{SA12}$ we also see that $h^{*}_{m,E,\Omega}$ is m-maximal in $\Omega\setminus E$ where $E$ is a compact set in $\Omega.$
According to Theorem 3.4 in \cite{SA12} (see also Theorem 2.24 in $\cite{Ch15}$), a Borel subset $E$ of $\Omega$ is $m$-polar if and only if $Cap_m(E)=0.$  

\n
A sequence $\{u_j\}$ of measurable functions on $\Om$ is said to converge to a function $u$ in $Cap_m$ if for every compact subset $K$ of $\Om$ and every $\ve>0$ we have
$$\lim_{j \to \infty }Cap_m \{z \in K: |u_j (z)-u(z)|>\ve \}=0.$$
\n {\bf 2.6.} Let $u\in SH_{m}(\Omega),$ and let ${\Omega_{j}}$ be a fundamental sequence of $\Omega,$ which means
$\Om_j$ is strictly $m-$ pseudoconvex, $\Omega_{j}\Subset\Omega_{j+1}$ and $\bigcup_{j=1}^{\infty}\Omega_{j}=\Omega.$ Set 
$$u^{j}(z)=\big(\sup\{\varphi(z):\varphi\in SH_{m}(\Omega), \varphi\leq u \,\,\text{on}\,\, \Omega_{j}^{c}\}\big)^{*},$$
where $\Omega_{j}^{c}$ denotes the complement of $\Omega_{j}$ on $\Omega.$\\
We can see that $u^{j}\in SH_{m}(\Omega)$ and $u^{j}=u$ on $(\overline{\Omega_{j}})^{c}.$ From definition of $u^{j}$ we see that $\{u^{j}\}$ is an increasing sequence and therefore $\lim\limits_{j\to\infty}u^{j}$ exists everywhere except on an $m-$ polar subset on $\Omega.$ Hence, the function $\tilde{u}$ defined by $\tilde{u}=\big(\lim\limits_{j\to\infty}u^{j}\big)^{*}$
is an $m-$ subharmonic function on $\Omega.$ Obviously, we have $\tilde{u}\geq u.$ Moreover, if $u\in \mathcal{E}_{m}(\Omega)$ then $\tilde{u}\in  \mathcal{E}_{m}(\Omega)$ and $\tilde{u}\in MSH_{m}(\Omega).$
Set
$$\mathcal{N}_{m}=\mathcal{N}_{m}(\Omega)=\{u\in\mathcal{E}_{m}(\Omega): \tilde{u} =0.\}$$
We have the following inclusion
$$\mathcal{F}_{m}(\Omega)\subset \mathcal{N}_{m}(\Omega) \subset \mathcal{E}_{m}(\Omega).$$

\n {\bf 2.7.} We recall some results on weighted $m-$ energy classes in $\cite{DT23}.$ Let $\chi: \mathbb{R}^{-}\longrightarrow \mathbb{R}^{-},$ where $\mathbb{R}^{-}:= (-\infty, 0]$,
be an increasing function. We put
$$\mathcal{E}_{m,\chi}(\Omega)=\{u\in SH_{m}(\Omega:\exists(u_{j})\in\mathcal{E}_{m}^{0}(\Omega), u_{j}\searrow u, \sup_{j}\int_{\Omega} -\chi(u_{j})H_{m}(u_{j})<+\infty\}.$$
Note that the weighted $m-$ energy class generalizes the energy Cegrell class $\mathcal{F}_{m,p},\mathcal{F}_{m}.$
\begin{itemize}
	\item When $\chi\equiv -1,$ then $\mathcal{E}_{m,\chi}(\Omega)$ is the class $\mathcal{F}_{m}(\Omega).$\\
	\item When $\chi(t) =-(-t)^{p},$ then $\mathcal{E}_{m,\chi}(\Omega)$ is the class $\mathcal{E}_{m,p}(\Omega).$
\end{itemize}
According to Theorem 3.3 in $\cite{DT23}$, if $\chi\not\equiv 0$ then $\mathcal{E}_{m,\chi}(\Omega)\subset \mathcal{E}_{m}(\Omega)$ which means that the complex $m-$ Hessian operator is well - defined on class $\mathcal{E}_{m,\chi}(\Omega)$ and if $\chi(-t)<0$ for all $t>0$ then we have $\mathcal{E}_{m,\chi}(\Omega)\subset \mathcal{N}_{m}(\Omega).$  Moreover, by Theorem 3.6 in $\cite{DT23}$ we have if $u\in \mathcal{E}_{m,\chi}(\Omega)$ then the $m-$ weighted energy of $u$ is finite i.e., $\int_{\Omega}(-\chi\circ u)H_{m}(u)<+\infty.$

\n {\bf 2.8.} We will frequently use the following 
version of the comparison principle: 

\n 
{\it If $u, v \in \F_m (\Om), u \le v$ then  $\int\limits_{\Om} H_m (u) \ge \int\limits_{\Om} H_m (v).$}

\n 
The above result follows easily from the following monotonicity property which can be shown by applying repeatedly the integration by parts formula:

\n 
{\it If $u, v \in \F_m (\Om), u \le v$ and $\va \in \E^0_m$ then  $\int\limits_{\Om} (-\va )H_m (u) \ge \int\limits_{\Om} (-\va) H_m (v).$}

\n 
A related property is the following unicity of functions in the class $\F_m:$

\n 
{\it  If $u, v \in \F_m (\Om), u \le v$ and if $H_m (u)=H_m (v)$ then $u=v$ on $\Om.$}

\n 
{\bf 2.9.} We need the following useful approximation result in the class $\E_{\chi,m}.$
\begin{lemma} \label{echi}
	Let $\chi: \mathbb R^{-} \to \mathbb R^{-}$ be an increasing function such that $\chi (-t) < 0$ for every $t > 0$,
	and $u \in \mathcal \E_{\chi,m} (\Om) $. Then there exists a sequence $\{u_j\} \in \E^0_m (\Om)$
	such that:
	
	\n 
	(a) $u_j \downarrow u$ on $\Om;$
	
	\n 
	(b) $H_m (u_j)$ has compact support in $\Om$;
	
	\n 
	(c) $\lim\limits_{j\to \infty}\int\limits_{\Omega}-\chi(u_j)H_m(u_j)=
	\int\limits_{\Omega}-\chi(u)H_m(u).$
\end{lemma}
\begin{proof}
	The proof is a minor extension of Lemma 3.12 in \cite{DT23}. For the sake of completeness, we indicate some details. First, by Theorem 3.6 in \cite{DT23} we have $\int\limits_{\Omega} -\chi(u)H_m (u)<\infty.$
	Next, we fix an element $\rho \in \E^0_m (\Om)$
	and let $\Omega_j \uparrow \Omega$ be an increasing sequence of relatively compact subsets of $\Omega.$
	For $j \ge 1$ we set
	$\mu_j:= \ind_{\{u>j\rho\} \cap \Om_j} H_m (u).$ Then the measures $\mu_j$ have the following properties:
	
	\n 
	(i) $\mu_j$ has compact support in $\Om;$
	
	\n 
	(ii) $\mu_j \le \mu_{j+1} \le H_m (u);$ 
	
	\n 
	(iii) $\mu_j$ puts no mass on $m-$polar sets in $\Om;$ 
	
	\n 
	(iv) $\int\limits_\Om d\mu_j \le \int\limits_{\{u>j\rho\}} H_m(u) \le j^m \int\limits_{\Om} H_m (\rho) <\infty$ 
	(by Lemma 5.5 in $\cite{T19}$).
	
	\n 
	Thus using (ii), (iv) and the main Theorem in $\cite{Gasmi}$ we can find $u_j \in \mathcal{N}_m (\Om)$ such that $u_j \ge u$
	and $H_m (u_j)=\mu_j.$ In particular $H_m (u_j)$ has compact support in $\Om.$
	Again, by the Corollary 5.8 in $\cite{T19}$, (ii) and (iii) we get $u_j \geq u_{j+1}$ and 
	$u_j \ge j\rho$. Therefore $u_j \in \E^0_m (\Om).$ Set $v:= \lim\limits_{j \to \infty} u_j.$ Then $v \ge u$ and by 
	our construction of $u_j$ we have $H_m (v)=H_m (u).$ Hence $u=v$, and so $u_j \downarrow u$ on $\Om.$ 
	Finally, by Lebesgue monotone convergence Theorem we have
	$$\begin{aligned} 
	\lim\limits_{j\to \infty}\int\limits_{\Omega}-\chi(u_j)H_m(u_j)&=\lim\limits_{j\to \infty}\int\limits_{\Omega}-\chi(u_j)d\mu_j\\
	&=\lim\limits_{j\to \infty}\int\limits_{\Omega}-\chi(u_j)\ind_{\{u>j\rho\} \cap \Om_j} H_m (u)\\
	&= \int\limits_{\Omega}-\chi(u)H_m(u). 
	\end{aligned}$$
The proof is therefore completed.	
\end{proof}	
\section{Subextension in the class $\mathcal{E}_{m,\chi}$}
We start with a condition for membership of $\mathcal{F}_{m}(\Omega)$ in terms of the relative $m-$ Hessian capacity of sublevel sets.
\begin{lemma}\label{lm2.1}
Let $\varphi\in SH_{m}^{-}(\Omega).$ Assume that
	$$\varlimsup_{s\to 0}s^{m}Cap_{m}(\{z\in\Omega:\varphi\leq -s\},\Omega)<\infty.$$
Then $\varphi\in\mathcal{F}_{m}(\Omega).$	
\end{lemma}
\begin{proof}
By Theorem 3.1 in $\cite{Ch15}$, there exists a sequence of functions $\varphi_{j}\in \mathcal{E}_{m}^{0}(\Omega)\cap C(\Omega)$ such that $\varphi_{j}\searrow \varphi$ in $\Omega.$
	It remains to show that $\sup_{j}\int_{\Omega}H_{m}(\varphi_{j})<+\infty.$	
	Firstly, we will prove that
	\begin{equation} \label{as01}
		\int_{\{\varphi_{j}\leq -s\}}H_{m}(\varphi_{j})\leq s^{m}  Cap_{m}(\{\varphi_{j}\leq -s\},\Omega) \quad \forall s>0.
	\end{equation}
Indeed, let $0<s<s_{0}$ then we have $\max(\varphi_{j},-s_{0})=\varphi_{j}$ on the {\it open} set $\{\varphi_{j}>-s_{0}\}.$ Note that, $\{\varphi_{j}>-s_{0}\}$ is a neighborhood of the boundary of $\{\varphi_{j}\leq -s\}.$ Thus, we have
$$\int_{\{\varphi_{j}<-s\}}H_{m}(\max(\varphi_{j},-s_{0})) =\int_{\{\varphi_{j}<-s\}}H_{m}(\varphi_{j}) .$$
Hence we obtain
\begin{align*}
s_{0}^{-m}	\int_{\{\varphi_{j}<-s\}}H_{m}(\varphi_{j})&=s_{0}^{-m}	\int_{\{\varphi_{j}<-s\}}H_{m}(\max(\varphi_{j},-s_{0}))\\
&=\int_{\{\varphi_{j}<-s\}}H_{m}(\dfrac{\max(\varphi_{j},-s_{0})}{s_{0}})\\
&\leq Cap_{m}(\{\varphi_{j}\leq -s\},\Omega).
\end{align*}
Letting $s_{0}\searrow s $ we obtain (\ref{as01}).
Next we let $s\searrow 0$ in (\ref{as01}) to obtain 	
	$$\int_{\Omega}H_{m}(\varphi_{j})\leq \varlimsup_{s\to 0}s^{m}Cap_{m}(\{z\in\Omega:\varphi\leq -s\},\Omega)<\infty$$
	for all $j$. Thus $\varphi\in\mathcal{F}_{m}(\Omega)$ as desired.
\end{proof}
\n 
The following Lemma plays a crucial role in our proof. It is inspired by an analogous result for plurisubharmonic functions in Lemma 4.4 of $\cite{H08}$. Note, however that, here we do not assume $u$ is continuous and $H_m (u)$ may not be compactly supported on $\Om.$ We should say that some ideas of the proof are taken from Section 4 in $\cite{H08}$.
\begin{lemma}\label{lm3.2}
	Let $\Omega\subset\tilde{\Omega}\Subset \mathbb{C}^{n}$ be bounded $m-$ hyperconvex domains. Let $u\in\mathcal{E}^{0}_{m}(\Omega)$ be such that $K:= \overline{supp (H_m(u))}$ is included in $\tilde \Om.$ 
	Then the maximal subextension $\tilde{u}$ defined by 
	$$\tilde{u}:=\sup\{v\in SH^{-}_{m}(\tilde{\Omega}): v\leq u \,\,\text{on}\,\, \Omega\}$$
	has the following properties: 
	
	\n 
	(i)  $\tilde{u}\in \mathcal{E}^{0}_{m}(\tilde{\Omega});$
	
	\n 
	(ii) $H_{m}(\tilde{u})=0$ on $X:=(\{\tilde{u}<u\}\cap\Omega)\cup (\tilde{\Omega}\setminus\Omega);$ 
	
	\n 
	(iii) $H_{m}(\tilde{u})\leq \ind_{\Omega\cap\{\tilde{u}=u\}}H_{m}(u).$
\end{lemma}
\begin{proof}
	We split the proof into some steps.\\
{\bf Step 1.} We show $\tilde{u}\in \mathcal{E}^{0}_{m}(\tilde{\Omega}).$ For this we let 
	$U \Subset \tilde \Om$ be an open neighborhood of  the compact set $K$. Choose $\varphi\in\mathcal{E}^{0}_{m}(\tilde{\Omega})\cap C(\tilde{\Omega})$ such that $\varphi\leq \inf\limits_{U}\hat{u}.$ Set $\Omega^*:= \Omega \setminus K.$ Then we have $H_m (u)=0$ on $\Omega^*$
	and 
	$$\varliminf\limits_{z \to \partial \Om^*, z \in \Om^*} (u(z)-\va(z)) \ge 0.$$ By the comparison principle we obtain that $u \ge \va$ on $\Om^*.$ Hence $u \ge \va$ on $\Om.$ So $\tilde u \ge \va$ on $\tilde \Om.$
	Thus $\tilde{u}\in \mathcal{E}^{0}_{m}(\tilde{\Omega}).$
	
\n 
{\bf Step 2.}
We prove (ii) and (iii) under the additional assumption that $u$ is continuous on $\Om.$ Then by setting
	$$\hat u:=\begin{cases} u \ \text{on}\ \Om\\
		0 \ \text{on}\ \tilde \Om \setminus \Om
	\end{cases}$$
	we obtain a continuous non-positive function $\hat u$ on $\tilde \Om.$ It is easy to check that
	\begin{equation} \label{conx}
		\tilde u:= \sup\{v\in SH^{-}_{m}(\tilde{\Omega} ): v\leq \hat u \,\,\text{on}\,\, \tilde \Omega\}. 
		\end{equation}
Next, we will prove $H_{m}(\tilde{u})=0$ on the open set $X=\{\tilde{u}<\hat{u}\}.$ By Theorem 1.2 in $\cite{Bl1}$, it suffices to show $\tilde{u}$ is $m-$maximal on $X.$ For this, let $\psi\in SH_{m}(X)$ such that $\psi\leq \tilde{u}$ outside a compact subset $K$ of $X.$ Put
		$$\tilde \psi:=\begin{cases}
		\max(\psi,\tilde{u})& \text{on}\,\, X\\
		\tilde{u}&\text{on}\,\, \tilde{\Omega}\setminus X.
	\end{cases}$$
Then $\tilde \psi \in SH_m^{-} (\tilde \Om).$ Our trick is to modify $\tilde \psi$ to obtain a competitor of 
$\tilde u$ in (\ref{conx}). For this, we note that there is $\de>0$ so small such that
$$\de \psi +(1-\de) \tilde u \le \hat u \ \text{on}\  K.$$
To see this, it suffices to choose $\de \in (0,1)$ such that
$$\de \sup_K (-\hat u) \le (1-\de) \inf_K (\hat u-\tilde u).$$
Thus $\psi_\de:=(1-\delta)\tilde{u}+\delta\tilde \psi \leq \hat{u}$ on $\tilde{\Omega}.$ By (\ref{conx}) we get $\psi_\de \leq \tilde{u}$ on $\tilde{\Omega}.$ It follows that $\tilde \psi \leq \tilde{u}$ on $\tilde{\Omega}.$ So  $\psi\leq \tilde{u}$ on $X.$ Hence, $\tilde{u}$ is $m-$maximal on $X$ as desired.\\
To prove (iii), take a compact set $K\Subset \Omega\cap\{\tilde{u}=u\}.$ Then for $\varepsilon >0$ we have $K\Subset \{\tilde{u}+\varepsilon>u\}\cap\Omega.$ By Theorem 3.6 in $\cite{HP17}$ we have
\begin{align*}
\int_{K}H_{m}(\tilde{u})&=\int_{K}\ind_{\{\tilde{u}+\varepsilon>u\}}H_{m}(\tilde{u})\\
&=\int_{K}\ind_{\{\tilde{u}+\varepsilon>u\}}H_{m}(\max(\tilde{u}+\varepsilon,u))\\
&\leq \int_{K}H_{m}(\max(\tilde{u}+\varepsilon,u)).
\end{align*}
Observe that $H_{m}(\max(\tilde{u}+\varepsilon,u))$ is weakly convergent to $H_{m}(u)$ when $\varepsilon\to 0$ so
we have
$$\int_{K}H_{m}(\tilde{u}) \le \varlimsup_{\ve \to 0} \int_{K}H_{m}(\max(\tilde{u}+\varepsilon,u)) \le \int_K H_m(u).$$
Thus we obtain (iii).\\
{\bf Step 3.}
Now we want to remove the assumption on continuity of $u.$ To this end, by Theorem 3.1 in [Ch15]
there exists a sequence $\{u_j\} \in \mathcal{E}_{m}^{0}(\Omega)\cap C(\overline{\Omega})$ such that $u_{j}\searrow u$ in $\Omega.$ Let 
$$\tilde{u_{j}}=\sup\{\varphi\in SH_{m}^{-}(\tilde{\Omega}): \varphi\leq u_{j} \,\, \text{on}\,\,\Omega \}.$$ 
	Then we see that $\tilde{u_j} \ge \tilde u$ for all $j$, and $\tilde u_j \downarrow v \in SH_m^{-} (\tilde \Om)$ on $\tilde \Om.$ This implies $v \ge \tilde u$ on $\tilde \Om.$
	On the other hand, since 
	$$v \le \lim_{j \to \infty} u_j =u \ \text{on}\ \Om$$
	we infer that $v \le \tilde u$ on $\tilde \Om.$ Putting all this together we have shown that 
	$$\lim_{j \to \infty} \tilde u_j^{\downarrow}= \tilde u \ \text{on}\ \tilde \Om.$$
	Since $u_j$ is continuous, by Step 2 we have:
	
	\n 
	(a) $H_{m}(\tilde{u_j})=0$ on $(\{\tilde{u_j}<u_j\}\cap\Omega)\cup (\tilde{\Omega}\setminus\Omega);$ 
	
	\n 
	(b) $H_{m}(\tilde{u_j})\leq \ind_{\Omega\cap\{\tilde{u_j}=u_j\}}H_{m}(u_j).$

First, we show $ H_{m}(\tilde{u})=0 $ on $\{\tilde{u}<u\}\cap \Omega.$
Indeed, by Step 2 we have $H_{m}(\tilde{u_{j}})=0$ on $(\{\tilde{u_{j}}<u_{j}\}\cap\Omega)\cup (\tilde{\Omega}\setminus\Omega).$ So we have
$$(\max({u_{j},-s})-\max(\tilde{u}_{j},-s))H_{m}(\tilde{u_{j}})=0$$
for all $s>0.$ Letting $j\to +\infty$ we obtain
$$(\max({u,-s})-\max(\tilde{u},-s))H_{m}(\tilde{u})=0,$$
for all $s>0.$ Hence $H_{m}(\tilde{u})=0 $ on $\{\tilde{u}<-s<u\}\cap\Omega$ for all $s>0.$ 
Moreover, since $$\{\tilde{u}<u\}\cap\Omega=\bigcup_{s\in\mathbb{Q}^{+}}\{\tilde{u}<-s<u\}\cap\Omega$$ 
we infer $ H_{m}(\tilde{u})=0 $ on $\{\tilde{u}<u\}\cap \Omega.$ 
Next, we will show $H_m (\tilde u)=0$ on $\tilde \Om \cap \partial \Om$.
Fix $\ve>0$. Choose an open subset $V \Subset \Om$ such that $\int\limits_{\Om \setminus V} H_m (u)<\ve.$ Let $U \subset \tilde \Om \setminus V$ be an open neighborhood of $\tilde \Om \cap \partial \Om$.
Then from (a) and (b) and the weak convergence $H_m (u_j) \to H_m (u)$ we obtain
$$\begin{aligned} 
\int\limits_{\tilde \Om \cap \partial \Om}  H_m (\tilde u) \le  \int\limits_U  H_m (\tilde u) &\le \varlimsup_{j \to \infty} \int\limits_U  H_m (\tilde u_j)\\
& \le  \varlimsup_{j \to \infty} \int\limits_{\tilde \Om \setminus V} H_m (\tilde u_j)\\
&=\varlimsup_{j \to \infty}\int\limits_{\Om \setminus V} H_m (\tilde u_j)\\
 & \le \varlimsup_{j \to \infty} \int\limits_{\Om \setminus V} H_m (u_j)\\
&\le \varlimsup_{j \to \infty} \int\limits_{\Om} H_m (u_j)-\int\limits_V H_m(u)\\
&=\int\limits_{\Om \setminus V} H_m(u)<\ve.
\end{aligned}$$
By letting $\ve \to 0$ we have $H_m (\tilde u)=0$ on $\tilde \Om \cap \partial \Om$.
Finally, by repeating the same argument as we used in the last part of Step 2 we con clude that
$$H_{m}(\tilde{u}) \le \ind_{\Omega \cap \{\tilde u=u\}}H_{m}(u).$$ 
The proof is complete.
\end{proof}
\n
Now we are ready to provide a subextension result for the class  $\mathcal{E}_{m,\chi}.$ In the case  $\chi(t)= -(-t)^{p}$ we obtain an extension of Theorem 4.1 in $\cite{H08}$ from the class $\mathcal{E}_{p}$ to  $\mathcal{E}_{m,p}.$
Note that we are {\it not} assuming $\Om$ is relatively compact in $\tilde{\Om}.$ Thus our result is somewhat stronger than the main theorem in \cite{Be11}.
\begin{theorem}\label{th3.6}
	Let $\Omega\subset\tilde{\Omega}\Subset \mathbb{C}^{n}$ be bounded $m-$hyperconvex domains. Let $\chi:\mathbb{R}^{-}\to \mathbb{R}^{-}$ be an increasing function such that $\chi(-t)<0$ for all $t>0.$ 
Let $u\in\mathcal{E}_{m,\chi}(\Omega)$. Consider the maximal subextension $\tilde {u}$ defined by 
	$$\tilde{u}= \sup\{\varphi\in SH_{m}^{-}(\tilde{\Omega}): \varphi\leq u \,\, \text{on}\,\,\Omega \}.$$
	Then the following statements hold true:
	
	\n (i)
	$\tilde{u}\in \mathcal{E}_{m,\chi}(\tilde{\Omega}) $ and $\tilde{u}\leq u$ on $\Omega.$\\
	\n (ii) $H_m(\tilde{u})=0$ on $X:=(\{\tilde{u}<u\})\cap\Omega.$\\
	\n (iii)  $H_m(\tilde{u})\leq\ind_{\Omega\cap\{\tilde{u}=u\}}H_m(u)).$\\
	\n (iv) 
	$\int\limits_{\tilde{\Omega}}-\chi(\tilde{u})H_{m}(\tilde{u})\leq \int\limits_{\Omega}-\chi(u)H_{m}(u).$
\end{theorem}
\begin{proof}
First we use Lemma \ref{echi}  to get a sequence $u_j\in\mathcal{E}_{m}^{0}(\Omega)$ such that $H_m (u_j)$ has compact support in $\Om, u_j\searrow u$ and $$\lim\limits_{j\to \infty}\int_{\Omega}-\chi(u_j)H_m(u_j)=
\int_{\Omega}-\chi(u)H_m(u).$$	
Now for each $j\ge 1,$ let $\tilde{u_{j}}$ denote the maximal subextension of the function $u_{j}$ given by 
	Lemma $\ref{lm3.2}.$ Then $\tilde u_j$ satisfies the following properties:
	
\n	 
(a) $\tilde{u_{j}}\searrow \tilde{u}, \tilde{u_{j}}\in\mathcal{E}_{m}^{0}(\tilde{\Omega}), H_{m}(\tilde{u}_{j})=0$ on $(\{\tilde{u}_{j}<u_{j}\}\cap \Omega)\cup (\tilde{\Omega}\setminus\Omega);$ 

\n 
(b) $H_m(\tilde{u}_j)\leq \ind_{\Omega\cap\{\tilde{u}_j=u_j\}}H_m(u_j).$ 

It then follows from (a) and (b) that
\begin{align*}
	\int_{\tilde{\Omega}}-\chi(\tilde{u_{j}})H_{m}(\tilde{u_{j}})&=\int_{\{\tilde{u_{j}}=u_{j}\}\cap \Omega}-\chi(\tilde{u}_{j})H_{m}(\tilde{u_{j}})\\
	&\leq \int_{\{\tilde{u_{j}}=u_{j}\}\cap \Omega}-\chi(u_{j})H_{m}({u_{j}})\\
	&\leq \int_{\Omega} -\chi(u_{j})H_{m}({u_{j}}).
	\end{align*}
Here the second estimate follows from (a) and the fact that $\chi$ is increasing.	
So $\tilde{u}\in\mathcal{E}_{m,\chi}(\tilde{\Omega}).$ Obviously, $\tilde{u}\leq u$ on $\Omega.$ Thus, we get (i).\\
Next by Theorem 3.3 in $\cite{DT23}$ we have $\tilde{u}\in\mathcal{E}_m(\tilde{\Omega}).$ So $H_m(\tilde{u}_j)$ converges weakly to $H_m(\tilde{u})$ as $j\to \infty.$ Using the same argument as in Step 3 of Lemma 3.2 we get (ii) and (iii).
	Finally, by the above estimates we have 
$$\begin{aligned} 
\int\limits_{\tilde{\Omega}}-\chi(\tilde{u})H_{m}(\tilde{u}) &\le \lim_{j \to \infty} \int_{\tilde{\Omega}}-\chi(\tilde{u_{j}})H_{m}(\tilde{u_{j}})\\
&\le  \lim_{j \to \infty} \int_{\Omega} -\chi(u_{j})H_{m}({u_{j}})\\
&=\int_{\Omega}-\chi(u)H_m(u).
\end{aligned}$$
Thus we get the desired conclusion (iv).
\end{proof}

\section{Approximation  in $\F_m$ and $\mathcal{E}_{m,\chi}$}
The first main result of this section gives some characterizations for the approximation of elements in $\F_m (\Om)$ by
functions in the same class but  defined on larger domains to be possible. The main implication $(iv) \Rightarrow (i)$ seems to be new even in the plurisubharmonic case.
\begin{theorem}\label{th4.1}
Let $\Omega, \Omega_j$ be a bounded $m-$ hyperconvex domains in $\mathbb{C}^n.$
Assume that 
\begin{equation} \label{as1}
\Omega  \subset \Om_j, \ \forall j \ge 1.
\end{equation}
Consider the following assertions:

\n 
(i) For every  $u\in\mathcal{F}_{m}(\Omega),$ there is a sequence $u_{j}\in\mathcal{F}_{m}(\Omega_{j})$ such that:

\n 
$u_j \to u$ in $Cap_m$ on $\Omega, u_j \le u$ on $\Om,$ and $H_m (u_j) \le \ind_{\Om} H_m (u)$ for all $j.$

\n 
(ii) For every  $u\in\mathcal{F}_{m}(\Omega),$ there is a sequence $u_{j}\in\mathcal{F}_{m}(\Omega_{j})$ such that:

\n 
$u_j \uparrow u$ q.e. on $\Omega, u_j \ge \inf\limits_{\Om} u$ and $H_m (u_j) \le \ind_{\Om} H_m (u)$ for all $j.$

\n 
(iii)  $\Om_{j+1} \subset \Om_j$ for all $j \ge 1$ and
there exists a sequence $v_{j}\in\mathcal{F}_{m}(\Omega_{j})$ such that 
$(\varlimsup\limits_{j \to \infty} v_j)^*:=v \in \mathcal{F}_m (\Om).$

\n 
(iv) For every compact subset $K$ of $\Om$ we have
$$\lim_{j \to \infty}Cap_{m}(K,\Omega_{j})=Cap_{m}(K,\Omega).$$

\n 
(v) $\Om_{j+1} \subset \Om_j$ for all $j \ge 1$ and for every compact subset $K$ of $\Om$ we have
$$\lim_{j \to \infty}Cap_{m}(K,\Omega_{j})=Cap_{m}(K,\Omega).$$
\n 
(vi) For every open subset $U$ of $\Om$ we have
$$\lim_{j \to \infty}Cap_{m}(U,\Omega_{j})=Cap_{m}(U,\Omega).$$

\n 
Then the following implications hold true:

\n 
$(ii) \Rightarrow (iv), (ii) \Rightarrow (vi), (iii) \Rightarrow (iv), (iv) \Rightarrow (i), (v) \Rightarrow (ii)$.
\end{theorem}
\begin{remark}
We do not know if either of the implications $(i) \Rightarrow (iv)$ or $(iv) \Rightarrow (ii)$ is true.	
\end{remark}
\begin{proof} 
$(ii) \Rightarrow (iv)$ Fix a compact subset $K$ of $\Om.$ 
Then there is a sequence $u_{j}\in\mathcal{F}_{m}(\Omega_{j})$ such that:

\n 
$u_j \uparrow h^*_{K,m,\Om} \in \mathcal{F}_m (\Om)$ q.e. on $\Om, u_j \ge -1$ and $\text{supp}H_m (u_j) \subset K$ for all $j.$

\n 
Then we have $H_m (u_j) \to H_m (h^*_{K,m,\Om,})$ on $\Om.$ This implies that
$$\begin{aligned} 
\varlimsup_{j \to \infty} Cap_m (K, \Om_j) &\le Cap_m (K, \Om)\\
 &=\int_K H_m (h^*_{K,m,\Om})\\
&=\int_{\Om} H_m (h^*_{K,m,\Om})\\
&\le \varliminf_{j \to \infty} \int_{\Om} H_m (u_j)\\
&= \varliminf_{j \to \infty} \int_{K} H_m (u_j)\\
&\le \varliminf_{j \to \infty} Cap_m (K, \Om_j).
\end{aligned}$$
The proof is then complete.

\n 
$(ii) \Rightarrow (vi)$ First we suppose $Cap_m (U, \Om)<\infty.$
Then fix $\ve>0$.
Let $K$ be a compact subset of $U$ such that
$$Cap_m (K, \Om) \ge Cap_m (U, \Om) -\ve.$$
Since $h^*_{m,K,\Om} \in \F_m (\Om)$ we can find a sequence $h_j \in \F_m (\Om_j), h_j \ge -1, h_j \uparrow h^*_{m,K,\Om}$ q.e. on $\Om$.
Then since $H_m (h_j) \to H_m (h^*_{m,K,\Om})$ weakly on $\Om,$ we have
$$\begin{aligned}
\varliminf_{j \to \infty} Cap_m (U, \Om_j) &\ge \varliminf_{j \to \infty} \int_{U} H_m (h_j)\\
& \ge \int_U H_m (h^*_{m,K,\Om})\\
& \ge  \int_K H_m (h^*_{m,K,\Om})\\
&=Cap_m (K, \Om)\\
&\ge Cap_m (U, \Om) -\ve\\
&\ge \varlimsup_{j \to \infty} Cap_m (U, \Om_j)-\ve
\end{aligned}$$ 
By letting $\ve \to 0$ we obtain the desired conclusion. 

\n
Now we consider the case $Cap_m (U, \Om)=\infty.$ Then fix $N \ge 1.$ Let $K$ be a compact subset of $U$ such that
$Cap_m (K, \Om) \ge N.$ The proof is then finished by the same reasoning as the above case.

\n 
$(iii) \Rightarrow (iv)$ Fix a compact set $K \subset \Om.$ First we note that 
$$h^{*}_{m,K,\Omega}\in\mathcal{E}^{0}_{m}(\Omega) , -1 \leq h^{*}_{m,K,\Omega}\leq 0\ \text{and}\  \text{supp}H_{m}( h^{*}_{m,K,\Omega})\subset K.$$
 Let $h_{j}$ be the maximal subextension of $h^{*}_{m,K,\Omega} $ to $\Omega_{j}$ as in Lemma $\ref{lm3.2},$ i.e., 
$$ h_{j}(z)=\sup\{\varphi(z):\varphi\in SH_{m}^{-}(\Omega_{j}),\varphi|_{\Omega}\leq h^{*}_{m,K,\Omega} \}.$$ 
Then we have $-1 \le h_{j}\in SH_{m}^{-} (\Omega_j)$ and $H_{m}(h_{j} )\leq \ind_{\Omega}H_{m }(h^{*}_{m,K,\Omega})$ on $\Omega_{j}.$ 
Moreover, by the assumption that $v\in\mathcal{F}_m(\Omega)$, after multiplying $v$ with a fixed positive constant we may achieve that 
$v<-1$ on some fixed neighbourhood $U$ of $K$. Thus, by Hartogs' Lemma $v_j<-1$ on $U,$ and therefore
$v_{j}\leq h^{*}_{m,K,\Omega}$ on $\Omega.$ So $
v_{j}\leq h_{j}$ on $\Om_j.$
Set $$h(z):= \lim_{j \to \infty} h_j (z)\uparrow , \ z \in \Om.$$
 Then $h \ge v$ in $\Om,$ so $h \in \mathcal F_m (\Om).$
Note that, from $h_{j}\leq h^{*}_{m,K,\Omega} $ on $\Omega,$ we also have $h\leq h^{*}_{m,K,\Omega}$  on $\Om.$
Since
$H_m (h_j) \to H_m (h)$ weakly as $j \to \infty$ we obtain 
\begin{equation}  \label{condi2} 
H_m (h) \le H_{m }(h^{*}_{m,K,\Omega}) \ \text{on}\ \Om,
\end{equation} 
and the following chain of estimates
\begin{equation*} \int_{\Omega}H_{m}(h^{*}_{m,K,\Omega}) \le  \int_{\Omega}H_{m}(h)\leq  \varliminf_{j \to \infty} \int_{\Om} H_m (h_j)
\le \int_{\Omega}H_{m}(h^{*}_{m,K,\Omega}).
\end{equation*}
It follows that
\begin{equation*} \int_{\Omega}H_{m}(h^{*}_{m,K,\Omega})= \int_{\Omega}H_{m}(h).
\end{equation*}
Combining this with ($\ref{condi2}$) we get $H_m (h^{*}_{m,K,\Omega}) =H_m (h).$ By the unicity property of $\F_m (\Om)$ we obtain $h= h^{*}_{m,K,\Omega}$ on $\Om.$
Hence $H_{m}(h_{j}) \to H_{m}(h^{*}_{m,K,\Omega})$ weakly when $j\to \infty.$
Note that $ \text{supp}H_{m}(h_{j})\subset K$ and 
$\text{supp}H_{m}(h^{*}_{m,K,\Omega})\subset K$ so we have
$$Cap_m (K, \Om)=\int_\Om H_{m}(h^{*}_{m,K,\Omega}) \le  \varliminf\limits_{j\to\infty }\int_{\Om}H_{m}(h_{j})\leq \varliminf\limits_{j\to \infty } Cap_{m}(K,\Omega_{j}).$$
Combining this with the obvious estimate $$Cap_m (K, \Om) \ge  Cap_{m}(K,\Omega_{j}), \ \forall j \ge 1$$
we obtain the desired conclusion.

\n 
$(iv) \Rightarrow (i)$  Fix $u \in \F_m (\Om)$. 
Then for $j \ge 1,$ we let $u_{j}$ be the maximal subextension of $u$ to $\Omega_{j}$ defined as in Theorem $\ref{th3.6}$ (with $\chi \equiv -1$)
$$u_{j}:=\sup\{\varphi\in SH_{m}^{-}(\Omega_{j}): \varphi|_{\Omega}\leq u\}.$$
We then have 
$$u_{j}\in\mathcal{F}_{m}(\Omega_{j}), H_{m}(u_{j})\leq \ind_{\Omega \cap \{u=u_j\}}H_{m}(u), u_j \le u \ \text{on}\ \Om.$$
We will now show that $u_j \to u$ in $Cap_m$ on $\Om.$
If this is not true, then by passing to a subsequence, we obtain $\ve>0$, a compact subset $L$ of $\Om$ 
	such that
	$$Cap_m \{z \in L: \vert u_j (z)-u(z)\vert>\ve\}>\ve, \ \forall j \ge j_0.$$
Hence
		\begin{equation} \label{cond1}
		Cap_m \{z \in L: u_j (z)<u(z)-\ve\}>\ve, \ \forall j \ge j_0.\end{equation}
Define the following $m-$ subharmonic functions on $\Om$
$$\tilde u_j:= (\sup_{k \ge j} u_k)^*,   v:= \lim_{j \to \infty} \tilde u_j.$$ 
Thus $v=(\varlimsup\limits_{j \to \infty} u_j)^*.$ 
A key step is to check $v \in \F_m (\Om).$ For this, we fix $s \in (0,1)$ and a compact subset $K$ of 
$\{z \in \Om: v(z) \le -s\}.$ Then by quasicontinuity of $\tilde u_j$ and $v,$ there exists a compact subset $K'$ of $K$
such that $\tilde u_j|_{K'}, v|_{K'}$ are continuous and $Cap_m (K \setminus K', \Om)<1.$
Then by Dini's Theorem $\tilde u_j$ converges uniformly to $v$ on $K'.$ So there exists $j_0$ such that
$$ K' \subset \{z \in \Om: \tilde u_j \le -s/2\}, \ \forall j \ge j_0.$$
Thus, we obtain the following chain of estimates
\begin{align*}
s^{m}Cap_{m}(K,\Omega)&\le s^{m}Cap_{m}(K',\Omega)+s^m\\
&=s^{m} [\lim\limits_{j\to\infty}Cap_{m}(K',\Omega_{j})+1]\\
&\leq s^{m}[\lim\limits_{j\to\infty}Cap_{m}(\{z\in\Omega_{j}: \tilde u_{j}(z)\leq -s/2\},\Omega_{j})+1]\\
&\leq s^{m}[\lim\limits_{j\to\infty}Cap_{m}(\{z\in\Omega_{j}: u_{j}(z)\leq -s/2\},\Omega_{j})+1]\\
&\leq 2^m\lim\limits_{j\to\infty}\int_{\Omega_{j}}H_{m}(u_{j})+s^m\\
&\leq 2^m\int_{\Omega}H_{m}(u)+1.
\end{align*}
Hence 
$$\varlimsup\limits_{s \to 0} s^{m}Cap_{m}(\{z \in \Om: v(z) \le -s\},\Omega) \leq 2^m\int_{\Omega}H_{m}(u)+1<\infty.$$
So by Lemma $\ref{lm2.1}$ we get that $v\in\mathcal{F}_{m}(\Omega).$ Observe that, since $u_{j}\leq u$ on $\Omega$ we have $v\leq u$ on $\Om.$
Now by Theorem 3.8 in $\cite{HP17}$ we obtain $H_{m}(\tilde u_{j}) \to H_{m}(v)$ weakly as $j\to \infty.$
Fix $\va \in \E^0_m,$ then from the above weak convergence and the monotonicity property we obtain 
$$\begin{aligned} 
\int\limits_{\Om} (-\va) H_m (u) &\le \int\limits_{\Om} (-\va) H_m (v)\\
& \le \varliminf\limits_{j \to \infty} \int\limits_{\Om} (-\va) H_m (\tilde u_j)\\
& \le \varliminf\limits_{j \to \infty} \int\limits_{\Om} (-\va) H_m (u_j)\\
&\le \int\limits_{\Om} (-\va) H_m (u).
\end{aligned}$$ 
Hence $$\int\limits_{\Om} (-\va) H_m (u)=\int\limits_{\Om} (-\va) H_m (v), \ \forall \va \in \E^0_m.$$
It follows that $H_m (u)=H_m (v).$ By unicity property for $\F_m (\Om)$ we obtain 
$u =v$ on $\Om.$ This implies $u=\varlimsup\limits_{j \to \infty} u_j$ q.e. on $\Om.$
On the other hand, we observe that (\ref{cond1}) implies that 
the set $\{z: u(z)>\varlimsup\limits_{j \to \infty} u_j (z)\}$ is not $m-$polar. We reach a contradiction.
The desired conclusion now follows.

\n 
$(v) \Rightarrow (ii)$ The proof is similar to that of $(iv) \Rightarrow (i).$
Fix $u \in \F_m (\Om)$ and define
$$u_{j}:=\sup\{\varphi\in SH_{m}^{-}(\Omega_{j}): \varphi|_{\Omega}\leq u\}.$$
Then under the assumption that $\Om_j$ is decreasing we have $u_j$ is increasing. Moreover, 
$$u_j \in \F_m (\Om_j), u_j \le u, H_m (u_j) \le \ind_{\Om} H_m (u).$$
Set 
$v:= (\lim u_j)^*.$ By the same reasoning as in the implication $(iv) \Rightarrow (i)$
we can check that $v \in \F_m (\Om)$ and that $v \le u$ on $\Om.$ Now the key point is to show $H_m (u)=H_m (v).$ For this purpose, we employ an argument similar to $(iii) \Rightarrow (iv).$ By the weak convergence 
$H_m (u_j) \to H_m (v)$ we have $$H_m (v) \le \ind_{\Om} H_m (u).$$
On the other hand, we have the following estimates
$$\int_{\Om} H_m (u) \le \int_{\Om} H_m (v) \le \varliminf_{j \to \infty} \int_{\Om} H_m (u_j) \le \int_{\Om} H_m (u).$$
Here the first inequality follows from the comparison principle. Thus $H_m (u)=H_m (v).$ So the unicity principle
for $\F_m (\Om)$ yields that $u=v$ on $\Om.$ Thus $u_j \uparrow u$ q.e. on $\Om.$
\end{proof}	
The following Lemma describes a situation when (iii) of Theorem \ref{th4.1} holds.
\begin{lemma}\label{lm4.3}
	Let $\Omega \Subset \mathbb{C}^{n}$ be a strongly $m-$ hyperconvex domain and $\Omega_{j}$ be a decreasing sequence of $m-$ hyperconvex domain such that $\Omega=(\bigcap\Omega_{j})^{\circ}.$ Let $K$ be a compact subset of $\Omega.$ Then we have $Cap_{m}(K,\Omega_{j})\to Cap_{m}(K,\Omega),$ as $j\to \infty.$
\end{lemma}
\begin{proof} Set $u:= h^{*}_{m,K,\Omega} \in \F_m (\Om).$
	It is clear to see that $\{h^{*}_{m,K,\Omega_{j}}\}_{j}$ is an increasing sequence of function and $\lim\limits_{j\to \infty}h^{*}_{m,K,\Omega_{j}}\leq u.$
	Let $\rho$ be a negative $m-$subharmonic exhaustion function for $\Om$ which is defined on some domains $\Om' \Supset \Om.$
	Fix $\varepsilon>0$ and $0<c<\varepsilon.$ Set
	$$v:=
	\begin{cases}
	\max \{u-\varepsilon,\rho-c\}& \,\,\text{on}\,\,\Omega\\
	\rho-c&\,\,\text{on}\,\,\Omega^{'}\setminus\Omega.
	\end{cases}
	$$
	Take $j$ so large such that $\Omega_{j}\subset \{z\in\Omega^{'}:\rho<c\}.$ Then we have $v\in SH^{-}_{m}(\Omega_{j})$. By multiplying $\rho$ with a positive real number, we can assume that $K\subset \{\rho\leq -1\}$ then $v|_{K}\leq -1,$ so $u-\varepsilon\leq h^{*}_{m,K,\Omega_{j}}$ for all $j$ large enough.
	Letting $j\to \infty$ and $\varepsilon\to 0$, we obtain 
	$h^{*}_{m,K,\Omega}\leq \lim\limits_{j\to+\infty}h^{*}_{m,K,\Omega_{j}}\leq h^{*}_{m,K,\Omega}.$
	So we get $h^{*}_{m,K,\Omega_{j}} \uparrow u$ on $\Om.$
	The desired conclusion now follows from $(iii) \Rightarrow (iv)$ in Theorem \ref{th4.1}.
\end{proof}
Now we come to the next main result of this section.
\begin{theorem} \label{th4.4}
	Let $\Omega\Subset \mathbb{C}^{n}$ be an $m-$hyperconvex domain and $\Omega_{j}$ be a decreasing sequence of $m-$hyperconvex domains containing $\Omega.$ 
	Assume that 
	\begin{equation} \label{cond2}
	\lim\limits_{j\to\infty}Cap_{m}(U,\Omega_{j})=Cap_{m}(U,\Omega),
	\end{equation} for all open subset $U$ of $\Omega.$ 
	Let $\chi:\mathbb{R}^{-}\to \mathbb{R}^{-}$ be a increasing function such that $ \chi(-t)<0$ for all $t>0$ and $\chi\in C^{1}(\mathbb{R}^{-}).$ Then for every  $u\in\mathcal{E}_{m,\chi}(\Omega) ,$ there exists an increasing sequence of functions $u_{j}\in\mathcal{E}_{m,\chi}(\Omega_{j}) $ such that $\lim\limits_{j\to \infty}u_{j}=u$ on $\Omega$ q.e. on $\Om.$
	\end{theorem}
\begin{remark}
(i) The condition (\ref{cond2})	holds true if $\Om_j$ and $\Om$ satisfy the condition (iii) in Theorem \ref{th4.1}
or that of Lemma \ref{lm4.3}.

\n 
(ii) We do not know if Theorem \ref{th4.4} remains true without assuming $\chi \in C^{1}(\mathbb{R}^{-}).$
\end{remark}	
We require the following result which is analogous to Proposition 3.1 in \cite{CKZ}
\begin{lemma}\label{lmt18}
	For every $u\in\mathcal{N}_{m}(\Omega)$ and for any $s,t>0$ we have
	\begin{equation*}
	t^{m}Cap_{m}(\{u<-s-t\})\leq \int_{\{u<-s\}}H_{m}(u).
	\end{equation*}
\end{lemma}
Next,we recall Theorem 3.11 in $\cite{DT23}$
\begin{lemma}\label{lm4.10} Let  $\chi:\mathbb{R}^{-}\to \mathbb{R}^{-}$ be a increasing function such that $\chi\in C^{1}(\mathbb{R}^{-}).$
	Put $$\hat{\mathcal{E}}_{m,\chi}(\Omega)=\{\varphi\in SH^{-}_{m}(\Omega): \int_{0}^{\infty}t^{m}\chi^{'}(t)Cap_{m}(\{\varphi<-t\})dt<\infty\}.$$
	Then we have 
	$$ \hat{\mathcal{E}}_{m,\chi}(\Omega)\subset \mathcal{E}_{m,\chi}(\Omega).$$
\end{lemma}
\begin{proof} (of Theorem \ref{th4.4})
	Let $u\in\mathcal{E}_{m,\chi}(\Omega).$ For each $j\in\mathbb{N},$ let $u_{j}\in\mathcal{E}_{m,\chi}(\Omega_{j})$ be the maximal subextension of $u$ to $\Omega_{j}$ as in the Theorem $\ref{th3.6}$ which means that $$u_{j}:=\sup\{v\in SH_{m}^{-}(\Omega_{j}): v\leq u \,\,\text{on}\,\,\Omega\}.$$ Note that from $\chi(-t)<0$ for all $t>0$  we have $u_{j}\in\mathcal{E}_{m,\chi}(\Omega_{j})\subset \mathcal{N}_{m}(\Omega_{j}).$ 
	Since the sequence $\Omega_{j}$ is decreasing, we infer that $u_{j}$ is increasing on $\Om.$ Define 
	$\tilde{u}:=(u_{j} \uparrow )^{*}.$ Clearly, we have $\tilde{u}\leq u$ on $\Omega.$
	On the other hand, by the assumption ($\ref{cond2}$) we have 
	$$Cap_{m}(\{\tilde u<- 2t\},\Omega_{j}) \uparrow Cap_{m}(\{\tilde{u}<-2t\}, \Om), \ \forall t>0.$$
Thus, using Lebesgue monotone convergence Theorem
coupling with Lemma $\ref{lmt18}$ (for $t=s$)  we obtain
	\begin{align*}
	&\int_{0}^{\infty}t^{m}\chi^{'}(-t)Cap_{m}(\{\tilde{u}<-2t\},\Omega)dt\\
	&=\lim_{j\to \infty} \int_{0}^{\infty}t^{m}\chi^{'}(-t)Cap_{m}(\{\tilde{u}<-2t\},\Omega_{j})dt\\
	&\leq\lim_{j\to \infty}\int_{0}^{\infty}t^{m}\chi^{'}(-t)Cap_{m}(\{u_{j}<- 2t\},\Omega_{j})dt\\
	&\leq \lim_{j\to \infty}\int_{0}^{\infty}\chi^{'}(-t)\int_{\{u_{j}<-t\}}H_{m}(u_{j})dt\\
	&\leq \lim_{j\to \infty}\int_{\Omega_{j}}-\chi(u_{j})H_{m}(u_{j})\\
	&\leq \int_{\Omega}-\chi(u)H_{m}(u)\\
	&<\infty.
	\end{align*}
	So $\tilde{u}\in\hat{\mathcal{E}}_{m,\hat{\chi}}(\Omega)$ with $\hat{\chi}(t)=\chi(\dfrac{t}{2}).$ Note that from $\chi(-t)<0$ for all $t>0$ we also have $\hat{\chi}(-t)<0$ for all $t>0.$ Hence, Lemma $\ref{lm4.10}$  implies that $$\tilde{u}\in\hat{\mathcal{E}}_{m,\hat{\chi}}(\Omega)\subset\mathcal{E}_{m,\hat{\chi}}(\Omega)\subset \mathcal{N}_{m}(\Omega).$$ Moreover, $u\in\mathcal{N}_{m}(\Omega)$ and $\tilde{u}\leq u$ on $\Omega.$
Now we suppose that $\{\tilde u <u\} \ne \emptyset.$ Then we take a compact subset $K$ of $\{\tilde u <u\}$
with $Cap_m (K)>0.$
Set $h:= h^*_{K,m,\Om}.$ Then	
$h \in SH_m^{-} (\Om), -1 \le h<0$ and $\int_{\Omega}-hH_{m}(u)< \infty,$ we claim that
\begin{equation} \label{con3}
\int_{\Om} (-h) H_m (\tilde u) \le \int_{\Om} (-h) H_m ( u).\end{equation} 
Indeed, choose a sequence $u^{k}\in\mathcal{E}_{m}^{0}(\Omega)\cap C(\Omega), u^{k}\searrow u$  as in Theorem 3.1 in $\cite{Ch15}.$ For each $j\ge 1$, let $u^{k}_{j}$ be the maximal subextension of $u^{k}$  on $\Omega_{j}$ given by
$$ u^{k}_{j}:=\sup\{v\in SH_{m}^{-}(\Omega_{j}): v\leq u^{k}\,\,\text{on}\,\,\Omega\}.$$
We see that $u^{k}_{j}$ is decreasing sequence when $k$  increases to $\infty$ and $\lim\limits_{k\to \infty}u^{k}_{j}=u_{j}.$ By Lemma $\ref{lm3.2}$ we have $u^{k}_{j}\in\mathcal{E}_{m}^{0}(\Omega_{j})$ and $$H_{m}(u^{k}_{j})\leq \ind_{\Omega\cap\{u^{k}_{j}=u^{k}\}}H_{m}(u^{k}).$$
So using Theorem 3.8 in $\cite{HP17}$ we obtain
\begin{align*}
\int_{\Omega}-hH_{m}(\tilde{u})&=\lim\limits_{j\to \infty}\int_{\Omega}-hH_{m}(u_{j})\\
&=\lim\limits_{j\to \infty}\lim\limits_{k\to \infty}\int_{\Omega}-hH_{m}(u^{k}_{j})\\
&\leq \lim\limits_{k\to \infty}\int_{\Omega}-hH_{m}(u^{k})\\
&= \int_{\Omega}-hH_{m}(u).
\end{align*}
Thus, we have proved the claim (\ref{con3}). 	
Finally, by Lemma 2.9 in $\cite{Gasmi}$ we have 
	$$0 \le \fr1{m!} \int\limits_{\Omega} (u-\tilde u)^m H_m (h) \le \int\limits_{\Omega} 
	(-h) [H_m (\tilde u)-H_m (u)] \le 0.$$
	It follows that $$\int\limits_{\Omega} (u-\tilde u)^m H_m (h)=0.$$
	Hence $\int_K H_m (h)=0.$ This yields $Cap_m (K)=0.$ We obtain a contradiction to the choice of $K.$
The proof is therefore complete.
\end{proof}


\begin{thebibliography}{000000}

\bibitem[ACL18]{AL} P. \"Ahag, R. Cyz and L. Hed, {\it The geometry of $m-$hyperconvex domains,} The Journal of Geometric Analysis, {\bf 28}, (2018) 3196–3222.

\bibitem[BT82]{BT1} E. Bedford and B. A.Taylor, {\it A new capacity for plurisubharmonic functions}, Acta Math, {\bf 149} (1982), 1-40.

\bibitem[BT87]{BT87} E. Bedford and B. A.Taylor, {\it Fine topology, Silov boundary, and $(dd^c)^n$}, J. Funct. Anal. {\bf 72} (1987), 225-251. 
 
\bibitem[Be06]{Be06} S. Benelkourchi, {\it A note on the approximation of plurisubharmonic functions}, C. R.Acad. Sci. Paris, {\bf 342} (2006), 647-650.

\bibitem[Be11]{Be11} S. Benelkourchi, {\it Approximation of weakly singular of plurisubharmonic functions }, Int. J. Math, {\bf 22} (7) (2011), 937-946.

\bibitem[B{\l}05]{Bl1} Z. B{\l}ocki, {\it Weak solutions to the complex Hessian equation}, Ann. Inst. Fourier (Grenoble), {\bf 55} (2005), 1735-1756.
 

\bibitem[Ce98]{Ce98} U. Cegrell, {\it Pluricomplex energy}, Acta Math, {\bf 180} (1998), 187-217.

\bibitem[Ce04]{Ce04} U. Cegrell, {\it The general definition of the complex Monge-Amp\`ere operator}, Ann. Inst. Fourier (Grenoble), {\bf 54} (2004), 159-179.

\bibitem[CKZ05]{CKZ} U. Cegrell, S. Kolodziej and A. Zeriahi, {\it Subextension of plurisubharmonic functions
with weak singularities}, Math. Zeit., {\bf 250} (2005), 7-22.

\bibitem[CH08]{CH08} U. Cegrell and L. Hed, {\it Subextension and approximation of negative plurisubharmonic functions}, Michigan Math. J {\bf 56} (2008), no.3,  594-601.

\bibitem[CZ03]{CZ03} U. Cegrell and A. Zeriahi, {\it Subextension of plurisubharmonic functions with bounded complex Monge-Amp\`ere operator mass}, C. R.Acad. Sci. Paris Ser. I, {\bf 336} (2003), 305-308.

\bibitem[Ch12]{Ch12} L. H. Chinh, {\it On Cegrell's classes of $m$-subharmonic functions}, arXiv 1301.6502v1.

\bibitem[Ch15]{Ch15} L. H. Chinh, {\it A variational approach to complex Hessian equation in $\mathbb{C}^{n}$}, J. Math. Anal. Appl. , {\bf 431} (1) (2015), 228-259.

\bibitem[DBH14]{DBH14} LN. Q. Dieu, P. H. Bang and N. X. Hong, {\it Uniqueness properties of $m-$ subharmonic functions in Cegrell classes}, J. Math. Anal. Appl. , {\bf 420} (1) (2014), 669-683.



\bibitem [DT23]{DT23} D. T. Duong and N. V. Thien, {\it On the weighted $m-$ energy classes}, J. Math. Anal. Appl. , {\bf 519} (2) (2023), 126820.



\bibitem [Ga59]{Ga59} L. Garding, {\it An inequality  for hyperbolic polynomials }, J. Math and Mec., {\bf 8}  (1959), 957-965.

\bibitem[Ga21]{Gasmi} A.E Gasmy, {\it The Dirichlet problem for the complex Hessian operator in the class $\mathcal N_m (\Omega, f)$,} 
Math. Scand., {\bf 121} (2021), 287-316. 
\bibitem [H08]{H08} P. H. Hiep, {\it Pluripolar sets and the  subextension in Cegrell's classes}, Complex Var. Elliptic Equa., {\bf 53} (2008), 675-684.

\bibitem[HP17]{HP17} V. V. Hung and N. V. Phu, {\it Hessian measures on m- polar sets and applications to the complex Hessian equations}, Complex Var. Elliptic Equa., {\bf 62} (8) (2017), 1135-1164.


\bibitem[Kl91]{Kl} M. Klimek, Pluripotential Theory, The Clarendon Press Oxford University Press, New York, 1991, Oxford Science Publications.




\bibitem[SA12]{SA12} A. S. Sadullaev and B. I. Abdullaev, {\it Potential theory in the class of $m$-subharmonic functions}, Trudy Mathematicheskogo Instituta imeni V. A. Steklova, {\bf 279} (2012), 166-192.
\bibitem[T18]{T18} N. V. Thien, {\it A characterization of the Cegrell classes and generalized $m-$ capacities }, Ann. Pol. Math, {\bf 121} (1) (2018), 33-43.


\bibitem[T19]{T19} N. V. Thien, {\it Maximal $m$- subharmonic functions and the Cegrell class $\mathcal{N}_{m}$ }, Indagationes Mathematicae, {\bf 30} (2019), Issue 4, 717-739.

\end{thebibliography}
\end{document}